\documentclass[1p]{amsart}
\usepackage{amssymb}
\usepackage{amsmath}
\newtheorem{theorem}{Theorem}
\newtheorem{lemma}{Lemma}
\newtheorem{proposition}{Proposition}
\newtheorem{corollary}{Corollary}
\newtheorem{definition}[theorem]{Definition}

\newtheorem{example}[theorem]{Example}
\def\card#1{\vert#1\vert}
\def\cgd{\mathcal C \mathcal G \mathcal D}
\def\cgv{\mathcal C \mathcal G \mathcal V}
\def\cgp{\mathcal C \mathcal G \mathcal P}
\def\cgn{\mathcal C \mathcal G \mathcal N}
\def\cgsv{\mathcal C \mathcal G \mathcal S \mathcal V}
\def\cgsp{\mathcal C \mathcal G \mathcal S \mathcal P}
\def\cgvsv{\mathcal C \mathcal G \mathcal V \mathcal S \mathcal V}
\def\cgpsp{\mathcal C \mathcal G \mathcal P \mathcal S \mathcal P}
\def\cgvn{\mathcal C \mathcal G \mathcal V \mathcal N}
\def\cgpn{\mathcal C \mathcal G \mathcal P \mathcal N}

\begin{document}

\title{The golden number and Fibonacci sequences in the design of voting structures}
\author{Josep Freixas}
\address{Department of Applied Mathematics III and High Engineering School (Manresa Campus), Technical University of Catalonia, Spain. Tel.: +34-938777246, josep.freixas@upc.edu}
\author{Sascha Kurz}
\address{Department of Mathematics, University of Bayreuth, 95440 Bayreuth, Germany.\\ Tel.: +49-921-557353, Fax: +49-921-557352, sascha.kurz@uni-bayreuth.de }
\begin{abstract}
  Some distinguished types of voters, as vetoers, passers or nulls, as well as some others, play a significant
  role in voting systems because they are either the most powerful or the least powerful voters in the game
  independently of the measure used to evaluate power. In this paper we are concerned with the design of
  voting systems with at least one type of these extreme voters and with few types of equivalent voters, with this
  purpose in mind we enumerate  these   special classes of games and find out that its number always
  follows a Fibonacci sequence with smooth polynomial variations. As a consequence we find several families
  of games with the same asymptotic exponential behavior excepting of a multiplicative factor which is the golden
  number or its square. From a more general point of view, our studies are related with the design of voting
  structures with a predetermined importance ranking.
  
  \medskip
  
\end{abstract}

\maketitle

\section{Introduction}

Determining importance rankings is a significant issue in operational research. The structures which are the target of our study
completely rank the items (voters or components, see e.g.\ \cite{FP02ANOR}) from the most important to the least important according
to a well-known influence relation, so that we have a hierarchy for the items.  This total ranking for these structures implies also
the same ranking for the most well-known measures of importance~\cite{Fre10IJGT},~\cite{FMP12EJOR} and~\cite{DiMo02} so that it is
unchallengeable. 

In order to design structures or mechanisms for a given hierarchy we need to count all the possibilities available
for it. The main purpose  of this paper is enumerating these versatile structures commonly used in Operations Research. Indeed,
the study of ordinal preferences involves a variety of fields, including tournament theory, multiple criteria decision modeling
(MCDM), and, more recently, data envelopment analysis of qualitative data. As stated in the survey by Cook in~\cite{Coo06},
the notion of voter power or relative importance has been largely ignored in studies on ordinal ranking problems, although if a
tangible estimate of voter importance exists, then these voters can be treated like criteria in an MCDM problem. In fact, if a
common ranking exists for the most well-known power indices, this would definitively demonstrate a given importance ordering.
This approach is thus useful in designing structures ranking voters in voting institutions, workers in management enterprises
or device components. Examples in these different contexts can be found in:~\cite{ABCF09,JMK04,Lee02,ObIs03,SSCh09}.

Besides this more general motivation for our studies, the paper contributes to the classification of simple games or more generally
voting systems initiated in the classical monograph \cite{vNM44} 
by von Neumann and Morgenstern. Here we enumerate
some classes of complete simple games, i.e.\ special classes of voting systems in which each player casts a ``yes" or a ``no"
vote, and the outcome is a collective ``yes" or ``no" decision, with distinguished types of voters. We address our attention
to complete simple games with at least one of the six types of voters: dictators, veto players, passers, null players,
semi-veto players, or semi-passers, see Section~\ref{sec_technical_background} for the precise definitions. As far as we
know the two last types have not been considered before in voting literature.

\subsection{Related work}
One of the first and most prominent enumeration results on simple games, a super class of complete simple games, is \textit{May's
Theorem for Simple Games} \cite{May52} stating that the number of simple games with $n$ voters for simple majority decisions equals
one\footnote{May originally considered a slightly different setting of anonymous
voting systems for two alternatives that are \textit{neutral}, i.e.\ there is no built-in bias towards ``yes" or ``no" outcomes.}.
If only anonymous or symmetric (i.e.\ any pair of voters are substitutes) voters are considered for simple games with $n$ voters we 
get $$SG(n,1)=CG(n,1)=WG(n,1)=n.$$
Here $SG(n,t)$ denotes the number of simple games, $CG(n,t)$ the number of complete simple games, and $WG(n,t)$ the number of weighted
voting games with $n$ voters from $t$ different types of equivalent voters.

The number $CG(n,2)$ of complete simple games with $n$ voters belonging to exactly two types of voters were recently enumerated in
\cite{FMR12ANOR} and later on in \cite{KT12}, giving a simpler proof:
\begin{equation}
  \label{eq_cg_2}
  CG(n,2)=F(n+6)-(n^2 + 4n +8)\in \Theta\left(\left(\frac{1+\sqrt{5}}{2}\right)^n\right),
\end{equation}
where $F(n)$ are the Fibonacci numbers which constitute a well--known sequence of integer numbers defined by the following
recurrence relation: $F(0)=0$, $F(1)=1$, and $F(n)=F(n-1)+F(n-2)$ for all $n>1$.

The number of complete simple games with one shift-minimal winning coalition, see e.g.\ \cite{CaFr96,GHS12}, was
determined in \cite{FrPu98}: $\sum_{t=1}^n CG(n,t,1)=2^n-1$,
where $CG(n,t,r)$ denotes the number of complete simple games with $n$ voters, $t$ equivalent types of voters, and $r$ shift-minimal
winning coalitions. For complete simple games with two shift-minimal winning coalitions a more complicated enumeration formula was
determined in \cite{KT12}.
For given values of the parameters $t$ and $r$ it is possible to compute an exact enumeration formula for $CG(n,t,r)$
based on the parametric Barvinok algorithm and a tailored decomposition of a certain linear programming formulation for
complete simple games, see \cite{KT12}. We remark that the exact numbers of simple games are known up to $n=8$ voters and
the exact number of complete simple games or weighted voting games are known up to $n=9$ voters.

The structures under study, complete simple games, have interest in several different fields apart from voting although
we adopt in this paper the standard voting background. Fields for which these structures are of interest are: circuits,
clusters, threshold logic, cryptography, reliability or neural networks among others, see e.g.\ \cite{TaZw99} for an
overview. Recently simple games were studied using binary decision diagrams, see e.g.\ \cite{BBRdS11,Bo11}. 

Special types of voters in simple games were also considered in \cite{Pel81}. Complexity results for identifying some
of the proposed distinguished types of voters can be found in \cite{AzPhd}.

\subsection{Our contribution}
We establish bijections among several classes of complete simple games containing at least one of the mentioned
distinguished types of players, and obtain exact enumerations for these games with less than four types of
equivalent voters and for four types whenever null voters are present with either veto players or passers. 
While these enumerations are polynomial for only two types of voters, they follow a Fibonacci sequence modified
by a polynomial expression. So Fibonacci sequences and the golden number are in the core of these enumerations.
The obtained sequences in this paper have not yet appeared in the On-line Encyclopedia of Integer Sequences
(http://oeis.org/).

\subsection{Organization of the paper}
In Section~\ref{sec_technical_background} we briefly recall some necessary background on simple games, complete
simple games, i.e.\  the examined voting systems, and some types of players. To get the results in the present paper we
will make use of one additional previous result, besides the enumeration of $CG(n,2)$, namely  a parametrization
for complete simple games given in~\cite{CaFr96}, which we will restate in Subsection~\ref{subsec_characteristic invariants}.
In Section~\ref{sec_one_distinguished_type} we prove that four classes of complete simple games with voters being very
powerful or being nulls have all the same cardinality, additionally a closed formula on the number of voters is given
for the number of these games as long as the number of types of voters is less than~$4$. The main tool to that end
is Theorem~\ref{T:equivalences1} stating bijections between complete simple games containing at least one of the distinguished
types of voters. In Section~\ref{sec_two_distinguished_types} we consider the case of complete simple games containing
at least two of the distinguished types of voters. Finally, we end with a conclusion and some remarks on future research
in Section~\ref{sec_conclusion}.

\section{Technical Background}
\label{sec_technical_background}

In this section we briefly state the basic definitions of the examined voting systems.
For more background material we refer to \cite{TaZw99}.

One of the most general class of voting systems, in which a single alternative,
such as a bill or an amendment, is pitted against the status quo, is given by the
class of so-called simple games.


\begin{definition}\label{DSG}
A \emph{simple game} is a pair $(N, W)$ in which $N = \{1,2,\dots,n\}$
and $W$ is a collection of subsets of $N$ that
satisfies: $N \in W$, $\emptyset \notin W$ and, the monotonicity
property, if $S \in W$ and $S \subseteq T \subseteq N$ then $T \in W$.
\end{definition}

Any set of voters or players is called a \emph{coalition}, and the set $N$ is
called the \emph{grand coalition}. Members of $N$ are called
\textit{players} or \textit{voters}, and the subsets of $N$ that are
in $W$ are called \emph{winning coalitions}.
The subfamily of \emph{minimal} winning coalitions
$W^m = \{S\in W: T\subset S\Rightarrow T\notin W\}$
determines the game. The subsets of $N$
that are in $2^N\setminus{}W$ are called \emph{losing coalitions}.
Real--world examples of simple games are given in~\cite{Lee02,Tay95,TaZw99}
among others.

\begin{example}
  \label{ex_1}
  For $n=3$ we choose $W=\Bigl\{\{1,2\},\{1,3\},\{1,2,3\}\Bigr\}$, i.e.\ the set of winning coalitions, and $N=\{1,2,3\}$
  so that the set of minimal winning coalitions is given by $W^m=\Bigl\{\{1,2\},\{1,3\}\Bigr\}$.
\end{example}


The monotonicity property of simple games is some kind of a universal
base assumption for voting systems. Most of the systems used in practice
satisfy additional requirements. One common idea is the concept of influence,
i.e.\ that a particular voting system may give one voter more influence
than another. The so-called ``desirability" relation defined on the set of
voters represents a way to make this precise. Isbell already used it in
\cite{Isb58}.


\begin{definition}\label{Ddes}
Let $(N,W)$ be a simple game.
\begin{itemize}
\item[$(i)$] Player $i$ is \emph{at least as desirable} as $j$
($i \succsim  j$, in short) in $(N,W)$ if
$$ S\cup\{j\}\in W\ \Rightarrow   \ S\cup\{i\}\in W,\qquad
   \text{for all} \; \; S\subseteq N \setminus \{i,j\}.$$
\item[$(ii)$] Players $i$ and $j$ are \emph{equally desirable}
($i \approx j$, in short) in $(N,W)$ if
$$ S\cup\{i\}\in W \Leftrightarrow\ S\cup\{j\}\in W,\qquad
   \text{for all} \; \,  S \subseteq N\setminus \{i,j\}.$$
\item[$(iii)$] Player $i$ is \emph{strictly more desirable}
than player $j$ ($i \succ j$, in short) in $(N,W)$ if $i$ is at least
as desirable as $j$, but $i$ and $j$ are not equally desirable.
\end{itemize}
\end{definition}

In Example~\ref{ex_1} players~$2$ and $3$ are equally desirable, while
player~$1$ is strictly more desirable than player~$2$ and $3$. The 
$\approx$-relation partitions the set of voters into equivalence classes
$N_i$. We say that two voters in the same equivalence class have the
same influence and belong to the same \emph{type} (of voters). We also
speak of the number of types of voters meaning the number of
equivalence classes. W.l.o.g.\ we assume 
$1\succsim\dots\succsim n$ in the following, i.e.\ completeness of the
desirability relation on $(N,W)$.

\begin{definition}
\label{Dcomp} 
A simple game $(N,W)$ is \emph{complete} or \emph{linear} if the
desirability relation is a complete preordering.
\end{definition}

A system to amend the Canadian Constitution used in the sixties and 
studied in Kilgour~\cite{Kil83} is an example of a complete simple game
with two types of voters. While most of the common voting systems
consist of many types of voters, examples with only few types of
voters are not that artificial as one might think at the first
moment. As the structural properties of voting systems with few
types of voters seem to be more accessible, they are also a good
starting point for theoretical considerations. The enumeration of such
voting structures is indeed the target of our study. In the field of
Boolean algebra, complete simple games correspond to $2$-monotonic
positive Boolean functions, which were already considered in \cite{Hu65}.
The problem of identifying this type of functions by using polynomial-time
recognition has been treated in \cite{BHIK91,BHIK97}.


We introduce some distinguished types of voters for simple games.
\begin{definition}
\label{D:4types}
Let $(N,W)$ be a simple game.
\begin{itemize}
\item [$i)$] $i \in N$ is a \emph{dictator} in $(N,W)$ if $W^m = \bigl\{\{ i \}\bigr\}$,
\item [$ii)$] $i \in N$ has \emph{veto} in $(N,W)$ if $i \in S$ for all $S \in W$,
\item [$iii)$] $i \in N$ is a \emph{passer} in $(N,W)$ if $\{ i \} \in W$, 
\item [$iv)$] $i \in N$ is \emph{null} in $(N,W)$ if $i \notin S$ for all $ S \in W^m$.
\item [$v)$] $i \in N$ has \emph{semi-veto} in $(N,W)$ if $i \in S$ for all $S \in W$ 
      with the only exception of coalition $N \setminus \{i\}$ which is also a winning coalition.
\item [$vi)$] $i \in N$ is a \emph{semi-passer} in $(N,W)$ if $\{ i,j \} \in W$
      for all $j \in N \setminus \{i\}$, but $\{ i \} \notin W$.
\end{itemize}
\end{definition}
A typical example of an important complete voting system with vetoers
is the United Nations Security Council. The voters in this system are
the fifteen countries that make up the Security Council, five of which
are permanent members whereas the other ten are non-permanent members. 
Passage requires a total of at least nine of the fifteen possible votes,
subject to a veto due to a nay vote from any one of the five permanent
members. This model ignores abstention. For a treatment of this example
considering the possibility of abstention we refer the reader to~\cite{FrZw03}. 

In the former USSR the three top state officials, the President, the Prime minister,
and the Minister of Defence (Ustinov, Brezhnev, Kosygin), all had ``nuclear suitcases".
Any two of them could authorize a launch of a nuclear warhead. No one could do it
alone.\footnote{Example taken from the presentation ``Secret sharing schemes and
complete simple games" by Arkadii Slinko.} This situation can be modeled as a complete
simple game where the set of minimal winning coalitions is given by
$W^m=\bigl\{\{1,2\},\{1,3\},\{2,3\}\bigr\}$. Here all three players are both semi-vetoers
and semi-passers.

A real-world example containing a null voter is given by the early European Economic Community
between 1958 and 1973. The six founders were West Germany, France, Italy, the Netherlands,
Belgium, and Luxembourg. At that time the voting procedure could be represented as the
weighted voting game $[12;4,4,4,2,2,1]$, see the explanation for this notation at the
end of this section. One can easily check that Luxembourg is a null voter in this complete
simple game. See e.g.\ \cite{LeBrMoZa11} for analyses of all EU councils of ministers.

For further considerations on real-world voting systems we refer the
interested reader e.g.\ to \cite{FeMa98,Tay95,TaZw99}. 

If a game has a dictator, the dictator is the unique player with this status and the
remaining players are null voters. Being a dictator is the most radical form
of having veto and of being a passer. A player can have veto and be
a passer \emph{if and only if} the game is the dictatorship of this player. Thus, veto
and passers are pairwise incompatible in the same game if this is not the
dictatorship of a voter.  If some of these types of voters are present in the game,
then they form an equivalence class $N_1$ whose members dominate by the desirability
relation \emph{all} the other elements in $N$. On the other hand, it is obvious that
if $n=1$ there cannot be null voters; if $n>1$ and the game has null voters they
form an equivalence class $N_t$ whose members are dominated by the desirability
relation for any other player in $N$.

It is also obvious from the previous definition that veto and semi-veto (or
passer and semi-passer) voters can concur in the same game, while veto and
semi-passer or passer and semi-veto cannot concur in the same game (if we assume that
their roles are taken by different players). Semi-vetoers
belong to the same equivalence class $N_1$ which is the strongest one by the desirability
relation if the game has no veto players, while they belong to the second class $N_2$ if the
game has veto players, and similarly for semi-passers and passers. Finally,
the concurrence of semi-vetoers and nulls in the same game, and
similarly, the concurrence of semi-passers and nulls in the same game is not possible.

In summary, if a game has either a dictator, veto players or passers they
are the strongest players in the game and belong to the most powerful class; if a
game has null voters they are the weakest players in the game and belong to the least
powerful class of the game.\footnote{Almost all power indices respect the desirability
relation, e.g.\ the Shapley-Shubik, Banzhaf or Johnston indices~\cite{Fre10IJGT}, \cite{FMP12EJOR},
and ~\cite{DiMo02}. Thus belonging to the strongest class in the game means being
one of the most powerful players, and belonging to the weakest class in the game means being
one of the least powerful players in the game.} If a game has veto and semi-veto players or
passers and semi-passers, the semi-veto and semi-passers are the second
strongest players in the game, while in the absence of veto and passers they
are respectively the most strongest players in the game.

\begin{definition}
\label{Ddual}
The dual game $(N,W^{\ast})$ of a simple game $(N,W)$ is defined by
$W^{\ast} = \{S \subseteq N : N\setminus S \notin W \}$.
\end{definition}
Hence, to win in the dual game is to block in the original one ($S$ is \emph{blocking} in $(N,W)$ if
$N\setminus S \notin W$). It is easy to verify that: $(W^{\ast})^{\ast} = W$,
$\succsim^{\ast} \ = \ \succsim$ thus $\approx^{\ast} \ = \ \approx$, and a game is complete
\emph{if and only if} the dual is.

From the definition of dual game, it easily follows that:
\begin{enumerate}
\item If $(N,W)$ is the dictatorship of player $i$, then $(N,W)=(N,W^{\ast})$;
\item Voter $i\in N$ has veto in $(N,W)$  \emph{if and only if} voter $i\in N$ is a passer in $(N,W^{\ast})$;
\item Voter $i\in N$ is null in $(N,W)$  \emph{if and only if} voter $i\in N$ is null in $(N,W^{\ast})$.
\item Voter $i \in N$ has semi-veto in $(N,W)$ \emph{if and only if}  voter $i\in N$ is a semi-passer in $(N,W^{\ast})$.
\end{enumerate}

As a last refinement of voting systems we want to mention weighted voting games, 
which are a subclass of complete simple games. Here we assume that the voting
system can be represented by $n+1$ non-negative real numbers -- a quota $q$ and
weights $w_i$ for the players. As short-hand notation we use $[q;w_1,\dots,w_n]$.
A coalition $S$ is winning if and only if $\sum_{i\in S} w_i \ge q$. The 
voting system from Example~\ref{ex_1} is weighted and a possible realization
is given by $[3;2,1,1]$. 
We remark that we have $w_i>w_j$ if $i\succ j$, but not necessarily
$w_i=w_j$ if $i\approx j$ for two voters $i,j\in N$. Note that $[51;50,49,1]$ is
another realization for Example~\ref{ex_1} demonstrating the last property. 
However, as $w_i\ge w_j$ implies $i\succeq j$, independently of the representation
used for the weighted game, we can conclude that all weighted voting games are 
complete simple games. For more than $5$ voters one may easily find examples of
complete simple games not being weighted. Many of the real-world
voting systems are weighted voting games (or can at least be represented as the
intersection of a small number of weighted voting games). We remark that every
complete simple game is the intersection of a finite number, the minimal
number is called its \emph{dimension}, of weighted voting games, see \cite{TaZw99}
and \cite{FP08EJOR} for results on the dimension of complete simple games
with one shift-minimal winning coalition.

\subsection{Characteristic invariants of complete simple games}
\label{subsec_characteristic invariants}

In order to have a more compact notation for the set of winning coalitions
of a complete simple game we introduce models of coalitions, so-called
profiles. The notion is essential for our enumerations later on.

Now we are going to define the $\delta$-ordering, introduced in
\cite{CaFr96} for an arbitrary number of types of voters $t$, for this particular
case which will be useful to present Theorem \ref{T:CaFr96}
below.

\begin{definition}
\label{Dprofiles}
Let $n$ be the number of players of a complete simple game $(N,W)$
with a given number of equivalence  classes of voters $t$, $N_1 > 
\dots > N_t$, let $n_i = \card{N_i}$ for all $1 \leq i \leq t$,
then $(n_1, \dots, n_t) \in \mathbb{N}^t = \underbrace{\mathbb{N} \times \dots \times \mathbb{N}}_t$
with $\sum_{i=1}^t n_i=n$, the \emph{box (or hyper-rectangle) of $\prod_{i=1}^{t} (n_i+1)$ profiles}
for $(N,W)$ is:
$I_{n_1} \times \dots \times I_{n_t} = \{(m_1,\dots,m_t) \in (\mathbb{N} \cup \{ 0
\})^t \, : \, m_1 \leq n_1, \dots, m_t\leq
n_t \}.$
Let $\overline m = (m_1,\dots,m_t) \in I_{n_1} \times \dots \times I_{n_t}$, the \emph{$\delta$-ordering}
given by the comparison of partial sums in $I_{n_1} \times \dots \times I_{n_t}$ is:
$$ \overline p \, \delta \, \overline m \quad \emph{if and only if}
\quad \sum_{i=1}^k p_i \geq \sum_{i=1}^k m_i \; \, \text{for all} \; \, 1 \leq k \leq t.$$
\end{definition}

It is not difficult to check that the couple $(I_{n_1} \times \dots \times
I_{n_t}, \delta)$ is a distributive lattice that possesses a maximum
element $\overline n = (n_1, \dots, n_t)$, and a minimum element, $\overline 0 = (0,\dots,0)$. The profiles
in $I_{n_1} \times \dots \times I_{n_t}$ can be completely ordered by the
\emph{lexicographical ordering}; profile $\overline p$ is
lexicographically greater than $\overline m$ if  $\overline p \, \delta \, \overline m$.

In Example~\ref{ex_1} we have $t=2$, $N_1=\{1\}$, $N_2=\{2,3\}$, $n_1=1$, and $n_2=2$. With
this the model $(1,1)$ represents the coalitions $\{1,2\}$ and $\{1,3\}$. A more extensive
example is outlined later on.

The big utility of the profiles and the $\delta$-ordering is its compatibility
with the property of being a winning coalition. To be more precise, for a
given profile $\overline p$ all corresponding coalitions are either winning or
losing. So we also speak of winning or losing profiles. Moreover, for two
profiles with $\overline p \, \delta \, \overline m$ and $\overline{m}$
winning, we have that $\overline{p}$ must be winning. This permits to describe the
set of winning profiles using the set of winning profiles, which are minimal in the
$\delta$-ordering, only.

\begin{definition}
\label{Disom}
Two simple games $(N, W)$ and $(N', W')$ are said to be
\emph{isomorphic} if there is a bijective map $f:N\rightarrow N'$
such that $S \in W$ if and only if $f(S) \in W'$; $f$ is called an
isomorphism of simple games.
\end{definition}

In the following we restate a known parametrization result of complete
simple games up to isomorphism, which has three parts. The first part shows
how to associate a vector $\overline n$ and a matrix $\mathcal M$  to a
complete simple game $(N,W)$ and describes the restrictions that
these parameters need to fulfill. The second part establishes that
isomorphic complete simple games $(N,W)$ and $(N',W')$ correspond to
the same associated vector $\overline n$ and matrix $\mathcal M$
(uniqueness). The third part shows that a vector $\overline n$ and a
matrix $\mathcal M$ fulfilling the conditions in Part A correspond
to a complete simple game $(N,W)$ (existence). We also assume in what
follows that the rows (if $r\geq 2$) of $\mathcal M$ are lexicographically
ordered by partial sums, i.e.\ if $p < q$, then there exists some
minimum element $k$ ($1 \leq k \leq t-1$) such that $m_{p,k} > m_{q,k}$
and $m_{p,i} = m_{q,i}$ for all $i<k$.

\begin{theorem} (Carreras and Freixas' Theorems 4.1 and 4.2  in \cite{CaFr96})
\label{T:CaFr96}
\begin{description}
       \item[Part A] Let $(N,W)$ be a complete simple
       game with $t$ nonempty equivalence classes of voters $N_1> \dots > N_t$, let $ (n_1,\dots,n_t)$ be
       the vector defined by their cardinalities. For each coalition $S \in
       W$ we consider the node or profile $(s_1,\dots,s_t) \in I_{n_1} \times \dots \times I_{n_t}$
       with components $s_k = \card{S \cap N_k}$ ($k=1,\dots,t$), and let $\mathcal
       M$ be the matrix
       $$\mathcal M=\left(%
\begin{array}{cccc}
  m_{1,1} & m_{1,2} & \hdots & m_{1,t} \\
\vdots & \vdots & \ddots & \vdots \\
  m_{r,1} & m_{r,2} & \hdots & m_{r,t}\\
\end{array}%
\right)$$ with $m_{1,1}>0$ whose $r$ rows are the nodes corresponding to
       winning profiles which are minimal in the
       $\delta$-ordering, let $\overline {m}_p = (m_{p,1}, \dots , m_{p,t})$ be the p$th$ row of $\mathcal M$. Matrix $\mathcal M$
       satisfies the four conditions below:
\begin{enumerate}
      \item $n_k > 0$ for all $k=1, \dots, t$;
      \item $\overline 0 \leq \overline{m}_p \leq \overline n$ for $p =1, \dots, r$;
      \item $\overline m_p$ and $\overline m_q$ are not $\delta$-comparable for all $p \neq q$; and
      \item if $t>1$ then for every $k<t$ there exists some $p$ such that\footnote{The
                             lexicographic ordering chosen
                             guarantees uniqueness under
                             permutation of rows.
                             This lexicographic ordering is a plausible election which
                             could be replaced for other alternative criteria.}
      $$ m_{p,k}>0, \; \; m_{p,k+1} < n_{k+1}. $$
\end{enumerate}
     \item[Part B] (Uniqueness) Two complete simple games with two classes of voters
     $(N, W)$ and $(N', \mathcal
     W')$ are isomorphic if and only if $(n_1,n_2) = (n'_1,n'_2)$
     and $\mathcal M = \mathcal M'$.
     \item[Part C] (Existence) Given a vector $(n_1,n_2)$ and a
     matrix $\mathcal M$ satisfying the conditions of Part A, there
     exists a complete simple game $(N,W)$ with two classes of voters associated to
     vector $(n_1,n_2)$ and matrix $\mathcal M$.
\end{description}
\end{theorem}

To illustrate how $(N,W)$ is obtained from the pair $((n_1,n_2),
\mathcal M)$ we consider the following example:

\begin{example}
  \label{ex_2} Let $(n_1,n_2)=(2,3)$ and $\mathcal M = \left(%
\begin{array}{cc}
  2 & 0 \\
  0 & 3 \\
\end{array}%
\right).$  
\end{example}
The set of winning profiles in $I_{n_1} \times I_{n_2}$ is
$$\{(0,3), (1,2), (1,3), (2,0), (2,1), (2,2), (2,3) \}$$ because each of
these profiles either $\delta$-dominates $(2,0)$ or $(0,3)$; and the
set of minimal winning profiles is $\{ (0,3), (1,2), (2,0) \}$
because the other winning profiles can be obtained from one of these
three profiles by simply adding some elements in one or two components.
Those winning profiles, which are minimal in the $\delta$-ordering, are
called shift-minimal winning coalitions. Here the set of shift-minimal 
winning coalitions is given by $\{ (0,3), (2,0) \}$, e.g.\ $(1,2)$
$\delta$-dominates $(0,3)$ since the later one arises from a \textit{right-shift}
of one voter. If we take $N_1=\{1,2\}$ and $N_2=\{3,4,5\}$ then $N =
N_1 \cup N_2$ and
$$W^m = \bigl\{ \{1,2\}, \{1,3,4\}, \{1,3,5\}, \{1,4,5\}, \{2,3,4\},
\{2,3,5\}, \{2,4,5\}, \{3,4,5\} \bigr\},$$ where the first coalition
corresponds to profile $(2,0)$, the last coalition to profile
$(0,3)$ and the remaining intermingle coalitions to profile $(1,2)$.
We remark that stated example is weighted with representation
$[6;3,3,2,2,2]$.

As abbreviation we use the notation $p\succeq q$ if profile $p$ $\delta$-dominates
profile $q$. If neither $p\succeq q$ nor $q\succeq p$, i.e.\ if profiles $p$ and
$q$ are not $\delta$-comparable, we write $p\bowtie q$.

Theorem~\ref{T:CaFr96} is a \emph{parametrization} theorem
because it allows one to enumerate all complete simple games up to
isomorphism by listing the possible values of certain invariants.
Based on this it is possible to compute an explicit enumeration formula
for the number of complete simple games with $n$ voters if the parameters
$t$ and $r$, i.e.\ the number of columns and rows of $\mathcal M$, are specified
(but arbitrary), see \cite{KT12}. This paper goes more deeply into the issue
of enumerations for special cases without assuming information on $r$.

An interest of complete simple games has emerged recently in the field of
Cryptography. Indeed, the access structure in a secret sharing (see e.g., Stinson~\cite{Sti92})
can also be modeled by a simple game. To this end Simmons~\cite{Sim90} introduced the
concept of a hierarchical access structure. Gvozdeva et al.~\cite{GHS12} study complete
simple games with one-shift minimal winning vector (i.e.\ with a unique row for matrix
$\mathcal M$) and observe that they are isomorphic concepts to conjunctive and disjunctive
(for the dual game) hierarchically access structures. Moreover, both conjunctive and
disjunctive hierarchically access structures have been proved to be ideal (Tassa~\cite{Tas07})
which means they can carry secure (i.e., not giving any information about the secret to
unauthorized coalitions). 

If a complete simple game has special types of voters then the characteristic invariants,
$(\overline n, \mathcal M)$, that define it (and therefore fulfill the properties
in Theorem~\ref{T:CaFr96}-(A)) have specific forms. In the next lemma, which will be
intensively used in the following section, we consider the specific form of the
characteristic invariants for complete simple games containing voters as defined in
Definition~\ref{D:4types}.


\begin{lemma} \label{L:4types}
Let $(\overline n, \mathcal M)$ be a complete simple game with at least two voters and vector $\overline n = (n_1,n_2,\dots,n_t)$ for some $t \geq 1$ and $\sum\limits_{i=1}^t n_i=n$:
\begin{itemize}
\item[$i)$] If the game has a dictator then $t=2$,
$ \overline n = (1,n-1) $ and
$\mathcal M = \left( 1 \; 0 \right).$
\item[$ii)$] If the game has $k$ veto players then $k=n_1$ and matrix $\mathcal M$ has the form: $\mathcal M = (n)$ if $t=1$, $\mathcal M = (n_1 \ a)$ with $a<n_2$ if $t=2$, and
\begin{equation} \label{Eq:veto}
\mathcal M = \left(
\begin{array}{cccc}
  n_1 & m_{1,2} & \dots & m_{1,t}\\
  n_1 & m_{2,2} & \dots & m_{2,t} \\
  \vdots & \vdots & \ddots & \vdots \\
  n_1 & m_{r,2} & \hdots & m_{r,t}
\end{array}
\right)
\end{equation}
with $m_{1,2}>0$ if $t>2$.
\item[$iii)$] If the game has $k$ passers then $k=n_1$ and matrix $\mathcal M$ has the form:
$ \mathcal M = (1)$ if $t=1$,
$$
\mathcal M = \left(
\begin{array}{cccc}
   1 & 0  \\
 0 & b
\end{array}
\right)
$$
with $n_2>1$ and $b>1$ if $t=2$, and
\begin{equation} \label{Eq:passer}
\mathcal M = \left(
\begin{array}{cccc}
   1 & 0 & \dots & 0 \\
 0 & m_{2,2} & \dots & m_{2,t} \\
  \vdots & \vdots & \ddots & \vdots \\
  0 & m_{r,2} & \hdots & m_{r,t}
\end{array}
\right)
\end{equation}
with $m_{2,2}>0$ if $t>2$.
\item[$iv)$] If the game has $k$ null voters then $t>1$, $k=n_t$ and matrix $\mathcal M$ has the form:
$$ \mathcal M = \left(
\begin{array}{ccccc}
  m_{1,1} & m_{1,2} & \dots & m_{1,t-1} & 0\\
  m_{2,1} & m_{2,2} & \dots & m_{2,t-1} & 0\\
  \vdots & \vdots & \ddots & \vdots &  \vdots \\
  m_{r,1} & m_{r,2} & \hdots & m_{r,t-1} & 0
\end{array}
\right)
$$
\item[$v)$] If the game has $k_1$ vetoers and $k_2$ semi-vetoers then $k_1=n_1$, $k_2=n_2$ and matrix $\mathcal M$ has the form: $ \mathcal M = ( n_1 \ n_2-1) $ if $t=2$, the form
$$ \mathcal M = \left(
\begin{array}{ccc}
  n_1 & n_2 & c \\
  n_1 & n_2-1 & n_3
\end{array}
\right)
$$
with $n_3>1$ and $c<n_3-1$ if $t=3$, and the form
$$ \mathcal M = \left(
\begin{array}{cccccc}
  n_1 & n_2 & m_{1,3} & \dots & m_{1,t-1} & m_{1,t} \\
  n_1 & n_2 & m_{2,3} & \dots & m_{2,t-1} & m_{2,t} \\
  \vdots & \vdots & \vdots & \ddots &  \vdots & \vdots \\
  n_1 & n_2 & m_{r-1,3} & \hdots & m_{r-1,t-1} & m_{r-1,t} \\
  n_1 & n_2-1 & n_3 & \hdots & n_{t-1} & n_t \\
\end{array}
\right)
$$
if $t>3$.

If the game has no vetoers but $k$ semi-vetoers then $k=n_1$ and matrix $\mathcal M$ has the form:
$\mathcal M = (n-1)$ if $t=1$,
$$ \mathcal M =
\left(
  \begin{array}{cc}
    n_1 & a \\
    n_1-1 & n_2
  \end{array}
\right)
$$
with $n_2>1$ and $a<n_2-1$ if $t=2$, and the form
$$ \mathcal M =
\left(
  \begin{array}{cccc}
    n_1 & m_{1,2} &\dots & m_{1,t} \\
    n_1 & m_{2,2} & \dots & m_{2,t} \\
    \vdots & \vdots  & \ddots & \vdots \\
    n_1 & m_{r-1,2} & \dots & m_{r-1,t} \\
     n_1-1 & n_2& \cdots & n_t
  \end{array}
\right)
$$
if $t>2$.

\item[$vi)$] If the game has $k_1$ passers and $k_2$ semi-passers then $k_1=n_1$, $k_2=n_2$ and matrix $\mathcal M$ has the form:
$$
  \mathcal{M}=
  \left(
    \begin{array}{cc}1&0\end{array}
  \right)
$$
if $t=2$ and $n_2=1$ or
$$ \mathcal M = \left(
\begin{array}{cc}
  1 & 0 \\
 0 & 2
\end{array}
\right)
$$
if $t=2$ and $n_2>1$,
$$ \mathcal M = \left(
\begin{array}{ccc}
  1 & 0 & 0\\
 0 & 1 & 1 \\
  0 & 0 & b \\
\end{array}
\right)
$$
with $b>2$ and  $n_3>2$ if $t=3$ or
$$ \mathcal M = \left(
\begin{array}{ccc}
  1 & 0 & 0\\
 0 & 1 & 1
\end{array}
\right)
$$
with $n_3\ge 2$ if $t=3$, and
$$ \mathcal M = \left(
\begin{array}{cccccc}
  1 & 0 & 0 & \dots  & 0& 0 \\
  0 & 1 & 0 & \dots & 0 &1 \\
  0 & 0 & m_{3,3} & \dots & m_{3,t-1} & m_{3,t} \\
  \vdots & \vdots & \vdots & \ddots &  \vdots & \vdots \\
  0 & 0 & m_{r,3} & \hdots & m_{r,t-1} & m_{r,t} \\
\end{array}
\right)
$$
if $t>3$.

If the game has no passers but $k$ semi-passers then $k=n_1$ and matrix $\mathcal M$ has the form:
$\mathcal M = (2)$ if $t=1$,
$$
\left(
  \begin{array}{cc}
    1 & 1 \\
    0 & c
  \end{array}
\right)
$$
with $n_2>2$ and $c\ge 3$ if $t=2$, and the form 
$$
\left(
  \begin{array}{ccccc}
    1 & 0 &\dots & 0 & 1 \\
    0 & m_{2,2} & \dots & m_{2,t-1} & m_{2,t} \\
    \vdots & \vdots  & \ddots & \vdots & \vdots \\
    0 & m_{r,2}& \cdots & m_{r,t-1} & m_{r,t}
  \end{array}
\right)
$$
if $t>2$.
\end{itemize}
\end{lemma}

\section{Enumerations for complete simple games with either the most powerful voters or the least powerful -- one distinguished type of voters}
\label{sec_one_distinguished_type}

In Definition~\ref{D:4types} we have exposed  six distinguished types of voters. In this section we will consider
complete simple games containing at least one of those six special types of voters. To this end we introduce some
notation in Table~\ref{table_six_cases_one_distinguished_type}. So $CGV(n)$ e.g.\ represents the number of complete
simple games having vetoers with $n$ voters. If the number of different types additionally is restricted
we denote the corresponding number by $CGV(n,t)$. If we want to address the respective set of objects instead 
of their number we use the corresponding curly literals, i.e.\ $\cgd(n)$ $\cgv(n)$, $\cgp(n)$, $\cgn(n)$ $\cgsv(n)$,
$\cgsp(n)$ and $\cgd(n,t)$ $\cgv(n,t)$, $\cgp(n,t)$, $\cgn(n,t)$ $\cgsv(n,t)$, $\cgsp(n,t)$.

\begin{table}
  \begin{center}
    \begin{tabular}{ccccc}
      \hline
      \textbf{case} & \textbf{contained type} & \textbf{\# games} & \textbf{\# games with $\mathbf{t}$ types} \\
      i   & dictator    & $CGD(n)$  & $CGD(n,t)$  \\
      ii  & vetoer      & $CGV(n)$  & $CGV(n,t)$  \\
      iii & passer      & $CGP(n)$  & $CGP(n,t)$  \\
      iv  & null        & $CGN(n)$  & $CGN(n,t)$  \\
      v   & semi-veto   & $CGSV(n)$ & $CGSV(n,t)$ \\
      vi  & semi-passer & $CGSP(n)$ & $CGSP(n,t)$ \\
      \hline
    \end{tabular}
    \caption{Number of complete simple games with one distinguished type of voters.}
    \label{table_six_cases_one_distinguished_type}
  \end{center}
\end{table}

The condition that a complete simple game contains a dictator is very restrictive. Due to our assumption on the ordering of the players, the
first player is a dictator and the set of minimal winning coalitions is given by $W^m=\bigl\{\{1\}\bigr\}$. All other players then
have to be null voters. Thus we have $t=2$, $n_1=1$, $n_2=n-1$ and $\mathcal{M}=\left(
\begin{array}{cc}1&0\end{array}\right)$ unless $n=1$.

\begin{lemma} 
\label{lemma_dictator_enumeration}
$CGD(n)=1$ for $n\ge 1$ and $GCD(n,t)=\left\{\begin{array}{rcl}1&\text{if}&t=1,\,n=1,\\1&\text{if}&t=2,\, n\ge 2,\\0&&\text{otherwise}.
\end{array}\right.$
\end{lemma}

The remaining five classes of complete simple games with one of the distinguished types of voters are pairwise in one--to--one
correspondence and therefore their cardinalities coincide.

\begin{theorem} \label{T:equivalences1}
For all positive integers $n$ and $t$ there is a bijection among the sets of complete simple games:
$$\cgv(n,t), \; \cgp(n,t), \; \cgsv(n,t), \; \cgsp(n,t)$$ and for all positive integers $n$ and $t>1$ there is
a bijection between any of these four classes of games and $\cgn(n,t)$.
\end{theorem}

\noindent \emph{Proof:} Assume $t=1$, then $\overline n = (n)$ for any complete simple game with $n$ voters. There is
only one complete simple game with veto $\mathcal M = (n)$, only one complete simple game with passers
$\mathcal M = (1)$, only one complete simple game with semi-veto players $\mathcal M = (n-1)$, and only one 
complete simple game with semi-passers $\mathcal M = (2)$.

Thus, the pairwise bijections between pairs of sets $\cgv(n,1)$, $\cgp(n,1)$, $\cgsv(n,1)$ and $\cgsp(n,1)$ are clear.

Assume from now on that $t>1$. We define a bijection $f$ from $\cgv(n,t)$ to $\cgn(n,t)$; a  bijection $g$ from
$\cgp(n,t)$ to $\cgn(n,t)$; a bijection $h$ from $\cgv(n,t)$ to $\cgsv(n,t)$ and a bijection $k$ from
$\cgp(n,t)$ to $\cgsp(n,t)$. Of course, any other bijection is a composition of some of these bijections or
their inverses.

\begin{enumerate}

\item \emph{Definition of a bijection $f$ between $\cgv(n,t)$ and $\cgn(n,t)$.}

Let $(\overline n, \mathcal M) \in \cgv(n,t)$ be the characteristic invariants of a complete simple game with
$t$~types of voters having a veto. Let $\overline n = (n_1,n_2, \dots,n_t)$ and $m_{i,j}$ be the components
of $\mathcal M$ as defined in Equation~\eqref{Eq:veto} with $1 \leq i \leq r$ and $1 \leq j \leq t$ and
$m_{\cdot,1} = n_1$ because the game has veto players.\footnote{A dot in the first subindex $m_{\cdot,j}$ means that $i$ is any value between $1$ and $r$. Analogously, a dot in the second subindex in $m_{i,\cdot}$ means that $j$ is any value between $1$ and $t$.}

In order to define the bijection $f$, we distinguish two separate cases:

\begin{itemize}
\item[$(i)$] The game has no null voters.

Then, $f$ sends vector $(n_1,n_2,n_3,\dots,n_t)$ to vector $(n_2,n_3,\dots,n_t,n_1)$ and matrix $\mathcal M$
to matrix $\mathcal M'$ where for all $1 \leq i \leq r$ and $1 \leq j \leq t-1$ the elements of  $\mathcal M'$ are
defined as $m'_{i,j} = m_{i,j+1}$, whereas $m'_{\cdot,t} = 0$. That is,

\begin{equation} \label{Eq:fromvetotonull}
\mathcal M' = \left(
\begin{array}{cccc}
  m_{1,2} & \dots & m_{1,t} & 0\\
  m_{2,2} & \dots & m_{2,t} & 0 \\
   \vdots & \ddots & \vdots & \vdots \\
   m_{r,2} & \hdots & m_{r,t} & 0
\end{array}
\right)
\end{equation}

\item[$(ii)$] The game has null voters.

Then, $f$ is the identity.
\end{itemize}

Now we need to check that:
\begin{enumerate}
\item The map $f$ is well--defined.
\item The map $f$ is injective.
\item The map $f$ is surjective.
\end{enumerate}

\vskip 0.3truecm

\noindent \emph{Well-defined}. We have to prove that  $f(\overline n, \mathcal M) \in \cgn(n,t)$. This is trivially
true if $(\overline n, \mathcal M)$ has null voters. Otherwise, let $f(\overline n, \mathcal M) = (\overline n',
\mathcal M')$ where $\mathcal M'$ is defined as in Equation~\eqref{Eq:fromvetotonull}, this pair fulfills the
conditions established in Theorem~\ref{T:CaFr96}-(A), $\overline n'$ has $t$ components and $m'_{\cdot,t} = 0$
guarantees the presence of null voters. Hence, $(\overline n',\mathcal M') \in  \cgn(n,t)$.

\vskip 0.3truecm

\noindent \emph{Injective}. Let $(\overline n, \mathcal M), (\overline m, \mathcal P) \in \cgv(n,t)$ with
$f(\overline n, \mathcal M) = f(\overline{m}, \mathcal P)$.  If $f(\overline n, \mathcal M)$ has null voters
then $(\overline n, \mathcal M) = f(\overline n, \mathcal M) = f(\overline m, \mathcal P) = (\overline m,
\mathcal P)$. Otherwise, let $f(\overline n, \mathcal M) = f(\overline m, \mathcal P) = (\overline h, \mathcal Q)$,
then $$\overline{h} = (h_1,h_2, \dots, h_{t-1},h_t) = (n_2,\dots,n_t,n_1) = (m_2,\dots,m_t,m_1),$$ thus
$n_j = m_j$ for all $1 \leq j \leq t$ and we have $\overline n = \overline m$; moreover
$$ q_{i,j} = m_{i,j+1} = p_{i,j+1} \; \, \text{for all} \; \, i, \; \; \text{and for all} \; \, j<t,$$
Thus $\mathcal M$ and $\mathcal P$ have the same dimensions and their components coincide, with the possible
exceptions of the components appearing in their respective first columns, but $(\overline n, \mathcal M),
(\overline m, \mathcal P) \in \cgv(n,t)$ guarantees that all of them are equal to $n_1$ and therefore
$\mathcal M = \mathcal P$.

\noindent \emph{Surjective}. Let $(\overline m, \mathcal P) \in \cgn(n,t)$. If it has veto players then
$(\overline m, \mathcal P) \in \cgv(n,t)$ and  $f(\overline m, \mathcal P) = (\overline m, \mathcal P)$.
Otherwise, consider  $(\overline n, \mathcal M)$ defined as follows: 
$$ n_1 = m_t, \; \, n_i = m_{i-1} \; \, \text{for} \; \, 1<i\leq t,$$
$$ m_{\cdot,1} = m_t, \; \, m_{i,j} = p_{i,j-1} \; \, \text{for} \; \, 1<j\leq t.$$
The pair $(\overline n, \mathcal M)$ fulfills the conditions of Theorem~\ref{T:CaFr96}-(A) and has veto
players, therefore $(\overline n, \mathcal M) \in \cgv(n,t)$; and $f(\overline n, \mathcal M) = (\overline m, \mathcal P)$.

\vskip 0.3truecm

\item \emph{Definition of a bijection $g$ from $\cgp(n,t)$ to $\cgn(n,t)$.}

Let $(\overline n, \mathcal M) \in \cgp(n,t)$ be the characteristic invariants of a complete simple game with $t$ types
having passers. Let $\overline n = (n_1,n_2, \dots,n_t)$ and $m_{i,j}$ be the components of $\mathcal M$ as defined 
in Equation~\eqref{Eq:passer} with $1 \leq i \leq r$ and $1 \leq j \leq t$ and $m_{1,1} = 1$, $m_{1,j} = 0$ if $j>1$,
$m_{i,1} = 0$ if $i>1$ because the game has passers.

In order to define the bijection $g$, we distinguish two separate cases:

\begin{itemize}
\item[$(i)$] The game has no null voters.

Then, $g$ sends vector $(n_1,n_2,n_3,\dots,n_t)$ to vector $(n_2,n_3,\dots,n_t,n_1)$ and matrix $\mathcal M$ to matrix
$\mathcal M'$ with $r-1$ rows, where for all $1 \leq i \leq r$ and $1 \leq j \leq t-1$ the elements of  $\mathcal M'$
are defined as $m'_{i,j} = m_{i+1,j+1}$, whereas $m'_{\cdot,t} = 0$. That is,

\begin{equation} \label{Eq:frompassertonull}
\mathcal M' = \left(
\begin{array}{cccc}
  m_{2,2} & \dots & m_{2,t} & 0\\
  m_{3,2} & \dots & m_{3,t} & 0 \\
   \vdots & \ddots & \vdots & \vdots \\
   m_{r,2} & \hdots & m_{r,t} & 0
\end{array}
\right)
\end{equation}

\item[$(ii)$] The game has null voters.

Then, $g$ is the identity.
\end{itemize}

\vskip 0.3truecm

\noindent \emph{Well-defined}. All the rows of $\mathcal M$ are not pairwise $\delta$-comparable because they fulfill
the conditions in Theorem~\ref{T:CaFr96}-(A), as the game has $n_1$ vetoers, the first common component $n_1$ of all
the rows has no effect on these comparisons. Thus, the rows of $\mathcal M'$ are not pairwise $\delta$-comparable and
therefore the pair $(\overline n', \mathcal M')$ fulfills the conditions in Theorem~\ref{T:CaFr96}-(A) and $g$ turns
the $n_1$ vetoers into $n_1$ nulls. Therefore, $(\overline n, \mathcal M') \in \cgn(n,t)$.

The proof that $g$ is injective and exhaustive follows the same guidelines as for $f$.

\vskip 0.3truecm

\item   \emph{Definition of a bijection $h$ from $\cgv(n,t)$ to $\cgsv(n,t)$.}

Let $(\overline n, \mathcal M) \in \cgv(n,t)$ be the characteristic invariants of a complete simple game with $t$ types of voters
having at least one veto player. Let $\overline n = (n_1,n_2, \dots,n_t)$ and $m_{i,j}$ be the components of $\mathcal M$ as defined
in Equation~\eqref{Eq:veto} with $1 \leq i \leq r$ and $1 \leq j \leq t$ and $m_{\cdot,1} = n_1$ because the game has veto players.

In order to define the bijection $h$, we distinguish two separate cases:

\begin{itemize}
\item[$(i)$] The game has no semi-vetoers.

Then, $h$ leaves invariant the vector $(n_1,\dots,n_t)$ and adds a last row in $\mathcal M$ which is  $(n_1-1, n_2, \dots,n_t)$
getting $\mathcal M'$, i.e.\ this addition converts the $n_1$ vetoers into $n_1$ semi-vetoers since $N \setminus \{i\}$ for all
$i \in N_1$ turns into a winning coalition.

\item[$(ii)$] The game has semi-vetoers.

Then, $h$ is the identity.
\end{itemize}

\vskip 0.3truecm

\noindent \emph{Well-defined}. Each row of $\mathcal M$, which is also a row of $\mathcal M'$, dominates the new last row of
$\mathcal M'$ because the  first component of the former rows are equal to $n_1$ which is greater than $n_1-1$. Moreover,
$\sum_{i=1}^t m_{k,i} < n-1$ for all $k=1, \dots, r$ since $\sum_{i=1}^t m_{k,i} = n-1$ would imply that coalition
$N \setminus \{i\}$ for all $i \in N_1$ would be winning in the original complete simple game $(\overline n, \mathcal M)$,
which is a contradiction with the non-existence of semi-vetoers. Thus, $h(\overline n, \mathcal M) = (\overline n, \mathcal M')$
fulfills the conditions of Theorem~\ref{T:CaFr96}-(A) and $h$ turns the $n_1$ vetoers into $n_1$ semi-vetoers. Therefore,
$(\overline n, \mathcal M') \in \cgsv(n,t)$.

The proof that $h$ is injective and exhaustive follows the same guidelines as for $f$.

\vskip 0.3truecm

\item   \emph{Definition of a bijection $k$ from $\cgp(n,t)$ to $\cgsp(n,t)$}

Let $(\overline n, \mathcal M) \in \cgp(n,t)$ be the characteristic invariants of a complete game with passers and with
$t$~types of voters. Let $\overline n = (n_1, \dots,n_t)$ and $m_{i,j}$ be the components of $\mathcal M$ as defined
in Equation~\eqref{Eq:passer} with $1 \leq i \leq r$ and $1 \leq j \leq t$.

In order to define the bijection $k$, we distinguish two separate cases:

\begin{itemize}
\item[$(i)$] The game has no semi-passers.

Then, $k$ leaves invariant vector $(n_1,\dots,n_t)$ and transform the first row in $\mathcal M$ which is
$(1, 0, \dots,0,0)$ into $(1,0,\dots,0,1)$, while the rest of rows keep invariant to get $\mathcal M'$.

\item[$(ii)$] The game has semi-passers.

Then, $k$ is the identity.
\end{itemize}

\vskip 0.3truecm

\noindent \emph{Well-defined}. Of course the first new row of $\mathcal M'$ dominates the other rows because
$(1,0,\dots,0,1)$ $\delta$-dominates $(1,0,\dots,0,0)$ which $\delta$-dominates the rest of the rows of $\mathcal M'$. Conversely,
as the game has no semi-passers it does not have any row of type $(0,1,0,\dots,0,1)$ and all of its rows $\delta$-dominate
$(0,1,0,\dots,0,1)$ and as $m_{i,1}=0$ for all $i>1$ it follows that each row $k=2, \dots r$ of $\mathcal M$ does not $\delta$-dominate
$(1,0,\dots,0,1)$, thus they are pairwise non-$\delta$-comparable. Thus, $k(\overline n, \mathcal M) = (\overline n, \mathcal M')$ fulfills
the conditions of Theorem~\ref{T:CaFr96}-(A) and $k$ turns the $n_1$ passers into $n_1$ semi-passers. Therefore,
$(\overline n, \mathcal M') \in \cgsp(n,t)$.

The proof that $k$ is injective and exhaustive follows the same guidelines as for $f$.
\end{enumerate}
\hfill $\square$

We remark that we had to exclude the set $\cgn(n)$ from the bijections in Theorem~\ref{T:equivalences1} for $t=1$
since there are no simple games with only null voters due to the definition of a simple game. 

\begin{corollary}
  \label{cor_1}
  For all positive integers $n$ and $t$ we have
  $$
    CGV(n,t)=CGP(n,t)=CGSV(n,t)=CGSP(n,t).
  $$
  For $t>1$ we additionally have $CGN(n,t)=CGV(n,t)$ and $CGN(n,1)=CGV(n,1)-1=0$.
\end{corollary}

Having the enumeration results for complete simple games with $n$ voters and at most $2$ types of voters at hand, we
can conclude enumeration formulas for complete simple games with $t<4$ types of voters, where at least one
of these equivalence classes corresponds to a distinguished type of voters. We state the results for $CGV(n,t)$ only.
Combining them with Corollary~\ref{cor_1} yields five simultaneous equivalent enumerations.
\begin{proposition} $\,$\\
\begin{enumerate}
\item $CGV(n,1)=1$,
\item $CGV(n,2)=\dfrac{n(n-1)}2$ (whenever $n \geq 2$),
\item $CGV(n,3)=F(n+7)-\dfrac 12(n^3+2n^2+13n+26)$ (whenever $n \geq 4$),
\newline where $F(n)$ are the Fibonacci numbers which form a sequence
of integer numbers defined by the following recurrence relation: $F(0)=0$,
$F(1)=1$, and $F(n)=F(n-1)+F(n-2)$ for all $n>1$.
\end{enumerate}
\end{proposition}

\noindent \emph{Proof:}
\begin{enumerate}
\item If $t=1$, both vector $\overline n$ and matrix $\mathcal M$ must coincide with the number of
      voters $n$. Therefore $CGV(n,1)=1$.
\item If $t=2$, let $\overline n=(n_1,n_2)$ for a given game. As the game must have veto players then all the
      entries in the first column must be $n_1$. As the rows of $\mathcal M$ must be non-comparable by the
      $\delta$-relation, then $\mathcal M$ only should have a row with the only requirement that $m_{1,2}$
      must be different from $n_2$, otherwise all players would be equivalent and therefore $t=1$, which would
      be a contradiction. Thus we have
      $$
        CGV(n,2)=\sum_{n_1=1}^{n-1}\sum_{m_{1,2}=0}^{\overset{=n_2-1}{\overbrace{n-n_1-1}}}1=\dfrac{n(n-1)}2.
      $$
\item If we remove the null voters from a complete simple game with $t$ types of voters (having at least one null voter), we
      obtain a complete simple game with $t-1$ types of voters without null voters. We remark that due to the definition of
      a simple game it cannot consist of null voters only. Thus we have
      $$
        CGV(n,3)=\sum_{i=1}^{n} CG(n-i,2)-CGV(n-i,2)=\sum_{i=1}^{n-2}CG(n-i,2)-CGV(n-i,2).
      $$
Inserting $CGV(n,2)=\dfrac{n(n-1)}2$ and $CG(n,2)=F(n+6)-(n^2 + 4n +8)$, see Equation~\ref{eq_cg_2}, yields

\begin{eqnarray*}
  CGV(n,3) &=& \sum_{k=2}^{n-1} \left[ F(k+6) - \dfrac{k(k-1)}2 - (k^2+4k+8) \right]\\
  &&CGV(n,3)= F(n+7)-\dfrac 12 (n^3+2n^2+13n+26)
\end{eqnarray*}\hfill $\square$
\end{enumerate}


Starting with $n=4$ the first ten numbers of the sequence $CGV(n,3)$ are:
$$ 2, 11, 37, 98, 225, 470, 919, 1713, 3082, 5400. $$

Asymptotically we have
$$ \lim_{n \rightarrow \infty} \dfrac{CGV(n,3)}{CG(n,2)} = \lim_{n \rightarrow \infty} \dfrac{F(n+7)+O(n^3)}{F(n+6)+O(n^2)} = \lim_{n \rightarrow \infty} \dfrac{F(n+7)}{F(n+6)} = \dfrac{1 + \sqrt{5}}{2}.$$
Thus, when $n$ is fix and high enough, the number of complete games with three types of voters having veto players (or semi-veto, passer, semi-passer or null voters) is almost equal to the number of complete games with only two types of voters multiplied by $\dfrac{1+\sqrt{5}}{2}$.

We remark that the bijections from Theorem~\ref{T:equivalences1} respect the property of being weighted, which is not too hard to prove but we
omit this part because it is not needed for the purposes of this paper. Moreover, the bijections keep games within the class
of so-called $\alpha$-roughly weighted games,
where the fraction of the weight of the heaviest losing coalition divided by the weight of lightest winning coalition is at most
$\alpha\ge 1$, see \cite{GHS12IJGT}, where the authors have introduced this as one of three hierarchies in order to classify simple games.
The special case $\alpha=1$ corresponds to the better known class of roughly-weighted games, see e.g.\ \cite{GoSl11}.

%

\section{Enumerations with two distinguished type of voters in a complete simple game}
\label{sec_two_distinguished_types}

Continuing the considerations from Section~\ref{sec_one_distinguished_type} we study complete simple games
containing at least two of the six distinguished types from Definition~\ref{D:4types}. As mentioned in the
beginning of Section~\ref{sec_one_distinguished_type} each complete simple game with a dictator contains
$n-1$ null voters. By $CGDN(n)$ we denote the number of complete simple games with $n$ voters containing
a dictator and at least one null and by $CGDN(n,t)$ we denote the number of these objects additionally
restricted to exactly $t$~types of voters. From Lemma~\ref{lemma_dictator_enumeration} we conclude:

\begin{lemma}
  \label{lemma_dictator_enumeration_2}
  For $n\ge 2$ we have $CGDN(n)=CGD(n)=1$ and $$CGDN(n,t)=CGD(n,t)=
  \left\{\begin{array}{rcl}1&\text{if}& t=2,\\0&&\text{otherwise}.\end{array}\right.$$ For $n\le 1$ all
  four counts are zero.
\end{lemma}

Each dictator has also a veto and is a passer. So for each subset $\mathcal{S}\subseteq\{V,P\}$ we have
$CGD\mathcal{S}N(n)=CGD\mathcal{S}(n)$ and $CGD\mathcal{S}N(n,t)=CGD\mathcal{S}(n,t)$, extending the
previously used notation for the number of complete simple games with presence of some distinguished
types of voters in a natural way.

If a complete simple game containing a dictator contains a semi-vetoer or a semi-passer then we
have $n=2$ and $t=2$. 

\begin{lemma}
  \label{lemma_dictator_enumeration_3}
  For two subsets $\mathcal{S}_1\subseteq\{V,P,N\}$ and
  $\mathcal{S}_2\subseteq\{SV,SP\}$ with $|\mathcal{S}_2|\ge 1$ we have
  $$
    CGD\mathcal{S}_1\mathcal{S}_2(n)=\left\{\begin{array}{rcl}1&\text{if}&n=2,\\0&&\text{otherwise}\end{array}\right.
    \quad\text{and}\quad 
    CGD\mathcal{S}_1\mathcal{S}_2(n,t)=\left\{\begin{array}{rcl}1&\text{if}&n=t=2,\\0&&\text{otherwise}.\end{array}\right.
  $$
\end{lemma}

If a complete simple game contains a passer and a veto then these roles must be taken by the same player,
which then is a unique dictator, and for each subset $\mathcal{S}\subseteq\{SV,SP,N\}$ we have
$CGVP\mathcal{S}(n)=CGD\mathcal{S}(n)$ and $CGVP\mathcal{S}(n,t)=CGD\mathcal{S}(n,t)$, compare the
enumeration formulas in Lemma~\ref{lemma_dictator_enumeration}, Lemma~\ref{lemma_dictator_enumeration_2}, and
Lemma~\ref{lemma_dictator_enumeration_3}.

Next we consider the simultaneous occurrence of a semi-vetoer and a semi-passer in a complete simple game.
By $CGSVSP(n)$ we denote the number of complete simple games with $n$ voters containing at least one semi-veto
and at least one semi-passer. By $CGSVSP(n,t)$ we denote the number of these objects additionally
restricted to exactly $t$~types of voters. Here, and in the following definitions, we permit the somewhat
artificial situation that the two types of distinguished voters could be taken by the same player.

\begin{lemma} For $n\ge 1$ we have $$CGSVSP(n)=\left\{\begin{array}{rcl}0&\text{if}&n=1,\\2&\text{if}&n=3,\\1&&\text{otherwise}\end{array}\right.$$
and $$CGSVSP(n,t)=\left\{\begin{array}{rcl}1&\text{if}&t=1,n=3,\\1&\text{if}&t=2,n\ge 2,\\0&&\text{otherwise}.\end{array}\right.$$
\end{lemma}

\noindent \emph{Proof:}
We assume that voter~$1$ is a semi-passer. From the definition of a semi-passer we conclude
$\{1\}\notin W$, $\{1,2\},\{1,3\},\dots,\{1,n\}\in W$. For $n\ge 4$ voter~$1$ is the unique voter being contained
in each winning coalition besides $N\backslash\{i\}$. Thus voter~$1$ also has to be the unique semi-veto and we
conclude $$W^m=\bigl\{\{1,2\},\{1,3\},\dots,\{1,n\},\{2,3,\dots,n\}\bigr\}.$$ Since this simple game has the
weighted representation $[n-1;n-2;\overset{n-1}{\overbrace{1,1,\dots,1}}]$ it is complete. For $n=1$ no simple
game containing a semi-veto and a semi-passer exists. If $n=2$ and voter~$1$ is a semi-passer, then we have
$\bigl\{\{1\}\bigr\}\notin W$ and $\bigl\{\{1,2\}\bigr\}\in W$. Since for each semi-veto player $i$ the coalition
$N\backslash\{i\}$ is winning, we conclude that voter $1$ is also a semi-veto player and $\bigl\{\{2\}\bigr\}$
is a winning coalition, so that we have $W^m=\bigl\{\{2\}\bigr\}$. Here player~$2$ is a dictator, in contrast
to our usual ordering of the players. We remark that player~$1$ is also a null voter and thus has three roles
in this example. For $n=3$ we have the simple games uniquely characterized by
$W^m=\bigl\{\{1,2\},\{1,3\},\{2,3\}\bigr\}$ and $W^m=\bigl\{\{1,2\},\{1,3\}\bigr\}$.
In the first example voter~$1$ is both the unique semi-passer and the unique semi-veto. In the second example
all three players are semi-passers and semi-vetoers. In the third example voter~$1$ is the unique semi-passer
and players $2$ and $3$ are the semi-vetoers. We remark that all simple games with at most three voters are complete
(and weighted).\hfill{$\square$}

Next we consider complete simple games containing vetoers and semi-passers. The corresponding counts are denoted
by $CGVSP(n)$ and $CGVSP(n,t)$.

\begin{lemma}
  \label{L:vsp}
  For $n\ge 1$ we have 
  $$
    CGVSP(n)=\left\{\begin{array}{rcl}0&\text{if}&n=1,\\2&\text{if}&n=2,\\1&&\text{otherwise}\end{array}\right.
    \quad\text{and}\quad
    CGVSP(n,t)=\left\{\begin{array}{rcl}1&\text{if}&t=1,\,n=2,\\1&\text{if}&t=2,\,n\ge 2,\\
    0&&\text{otherwise}.\end{array}\right.
  $$
\end{lemma}

\noindent \emph{Proof:} 
For $n=1$ no semi-passer is possible and for $n=2$ there are only three complete simple games. So we can easily check that
the complete simple game given by  $W^m=\bigl\{\{1\},\{2\}\bigr\}$ does contain neither a vetoer nor a semi-passer. In
the complete simple game given by $W^m=\bigl\{\{1\}\bigr\}$ voter~$1$ has a veto and voter~$2$ is a semi-passer, so that
we have $t=2$. In the complete simple game given by $W^m=\bigl\{\{1,2\}\bigr\}$ both players are semi-passers and vetoers,
so that we have $t=1$. 

For $n\ge 3$ we assume that voter~$1$ is a semi-passer so that $\{1,2\},\{1,3\},\dots,\{1,n\}\in W$ and $\{1\}\notin W$.
The only possible veto player is voter~$1$ and no further vetoers or semi-passers can be present. Thus voter~$1$ forms
its own equivalence class of voters $N_1$. Since player~$1$ has a veto we conclude 
$W^m=\bigl\{\{1,2\},\{1,3\},\dots,\{1,n\}\bigr\}$ and $t=2$.\hfill{$\square$}


Similarly we denote by $CGPSV(n)$ and $CGPSV(n,t)$ the counts for complete simple games containing passers and semi-vetoers.
\begin{lemma}
  \label{L:psv}
  For $n\ge 1$ we have 
  $$
    CGPSV(n)=\left\{\begin{array}{rcl}0&\text{if}&n=1,\\2&\text{if}&n=2,\\1&&\text{otherwise}\end{array}\right.
    \quad\text{and}\quad
    CGPSV(n,t)=\left\{\begin{array}{rcl}1&\text{if}&t=1,\,n=2,\\1&\text{if}&t=2,\,n\ge 2,\\
    0&&\text{otherwise}.\end{array}\right.
  $$
\end{lemma}

\noindent \emph{Proof:} 
For $n=1$ no semi-vetoer is possible and for $n=2$ there are only three complete simple games. Assuming w.l.o.g.\
that voter~$1$ is a passer the remaining possibilities are given by $W^m=\bigl\{\{1\},\{2\}\bigr\}$ and 
$W^m=\bigl\{\{1\}\bigr\}$. In the first mentioned complete simple game both players are passers and semi-vetoers so
that we have $t=1$. In the second case voter~$1$ is a passer and voter~$2$ a semi-vetoer so that we have $t=2$.

Due to the definition of a semi-vetoer for $n\ge 3$ the only possible semi-vetoer is player~$1$ and there can be
no other passers besides player~$1$. Thus we have $W^m=\bigl\{\{1\},\{2,\dots,n\}\bigr\}$.\hfill{$\square$}

Even more mapping to dual games is a bijection between the corresponding classes of complete simple games of
Lemma~\ref{L:vsp} and Lemma~\ref{L:psv}, since the pairs vetoers/passers and semi-vetoers/semi-passers are
interchanged, see Definition~\ref{Ddual} and the comment thereafter.

Due to the definition of a semi-passer (for $n\ge2$) a complete simple game cannot contain both a semi-passer and a null voter.
Similarly a complete simple game cannot contain both a semi-vetoer and a null voter.

The remaining pairs of distinguished types of voters are veto/semi-veto, passer/semi-passer, veto/null, and
passer/null. 
The respective counts are denoted by $CGVSV(n)$, $CGPSP(n)$, $CGVN(n)$, $CGPN(n)$, $CGVSV(n,t)$, $CGPSP(n,t)$,
$CGVN(n,t)$, and $CGPN(n,t)$. We denote the corresponding classes by $\cgvsv(n)$, $\cgpsp(n)$, $\cgvn(n)$,
$\cgpn(n)$, $\cgvsv(n,t)$, $\cgpsp(n,t)$, $\cgvn(n,t)$, and $\cgpn(n,t)$.

\begin{theorem} \label{T:bijection2} For all positive integers $n>1$ and $t>1$ there is a bijection between
$\cgvsv(n,t)$, $\cgpsp(n,t)$, $\cgvn(n,t)$ and $\cgpn(n,t)$.
\end{theorem}

\noindent \emph{Proof:} We first define a bijection $h'$ between $\cgvn(n,t)$ and $\cgpn(n,t)$ in the following way;
if $(\overline n, \mathcal M) \in \cgvn(n,t)$ 
$$h'(\overline n, \mathcal M) = (\overline n, \mathcal M^{\ast})$$
where $(\overline n, \mathcal M^{\ast})$ is the dual game of $(\overline n, \mathcal M)$.

Now we need to check that:
\begin{enumerate}
\item The map $h'$ is well--defined.
\item The map $h'$ is injective.
\item The map $h'$ is surjective.
\end{enumerate}

\noindent \emph{Well-defined}. We have to prove that  $h'(\overline n, \mathcal M)=(\overline n, \mathcal M^{\ast})
\in \cgpn(n,t)$. This is trivially true because the dual of a complete game is a complete game too, $i \in N$
is null in $(N,W)$ \emph{if and only if} $i \in N$ is null in $(N,W^{\ast})$, and
$i \in N$ has veto in $(N,W)$ \emph{if and only if} $i \in N$ is a passer in $(N,W^{\ast})$.

\vskip 0.3truecm

\noindent \emph{Injective}. Let $(\overline n, \mathcal M), (\overline m, \mathcal P) \in \cgvn(n,t)$ with
$h'(\overline n, \mathcal M) = h'(\overline{m}, \mathcal P)$.  Let
$h'(\overline n, \mathcal M) = h'(\overline{m}, \mathcal P) = (\overline{h}, \mathcal Q)$. Because of the definition of $h'$ we have:
$$\overline h = \overline n = \overline m \qquad \text{and} \qquad \mathcal Q = \mathcal M^{\ast} =  \mathcal P^{\ast}$$
and after applying duality to the last expression we obtain the desired equality, i.e.\
$\mathcal Q^{\ast} = ({\mathcal M}^{\ast})^{\ast} =  \mathcal M $ and
$\mathcal Q^{\ast} = ({\mathcal P}^{\ast})^{\ast} = \mathcal P$

\noindent \emph{Surjective}. Let $(\overline m, \mathcal P) \in \cgpn(n,t)$, then it is clear that
$(\overline m, {\mathcal P}^{\ast}) \in \cgvn(n,t)$ and $h'(\overline m, {\mathcal P}^{\ast}) =
(\overline m, \mathcal P)$.

\vskip 0.3truecm

Next we remark that the same mapping is also a bijection from $\cgvsv(n,t)$ to $\cgpsp(n,t)$. As before all three
conditions can be easily checked.

\vskip 0.3truecm

Finaly we define a bijection $h''$ between $\cgvsv(n,t)$ and $\cgvn(n,t)$ in the following way; if
$(\overline{n},\mathcal{M})\in \cgvsv(n,t)$, with $\overline{n}=(n_1,\dots,n_t)$ and $\mathcal{M}=(m_{i,j})_{1\le i\le r,1\le j\le t}$,
then $h''(\overline{n},\mathcal{M})=(\overline{n}',\mathcal{M}')$, where
\begin{eqnarray*}
  \overline{n}' &=& (n_1,n_3,\dots, n_t,n_2),\\
  \mathcal{M}'  &=& \begin{pmatrix}
                      m_{1,1} & m_{1,3} & \cdots & m_{1,t} & 0 \\
                      \vdots  & \vdots  & \ddots & \vdots    & \vdots \\
                      m_{r-1,1} & m_{r-1,3} & \cdots & m_{r-1,t} & 0
                    \end{pmatrix}.
\end{eqnarray*}
We note that we have $r\ge 2$. In other words, we have shifted the second equi\-va\-lence class of voters
to the $t$-th equi\-va\-lence class, while replacing the second
column $(m_{1,2},\dots,m_{r,2})^T$ by the all-zero vector and deleted the last row of $\mathcal{M}$.

Due to Lemma~\ref{L:4types}.v) the second column and the last row of $\mathcal{M}$ are uniquely characterized by the values of $n_1,\dots, n_t$.
So it remains to check that the resulting games are complete simple games with $t$ equivalence classes of voters, i.e.\ that they satisfy the
conditions from Theorem~\ref{T:CaFr96}, which can be done easily.\hfill $\square$

%

\begin{corollary}
  \label{cor_2}
  For all positive integers $n>1$ and $t>1$ we have
  $$
    CGVSV(n,t)=CGPSP(n,t)=CGVN(n,t)=CGPN(n,t).
  $$
\end{corollary}

\begin{proposition} \

\begin{enumerate}
\item $CGVN(n,2)=n-1$ (whenever $n \geq 2$),
\item $CGVN(n,3)=\dfrac{(n-1)(n-2)(n-3)}6$ (whenever $n \geq 4$),
\item $CGVN(n,4)=F(n+8)-\dfrac 16(n^4-2n^3+26n^2+47n+132)$ (whenever $n \geq 5$).
\end{enumerate}
\end{proposition}

\noindent \emph{Proof:}
\begin{enumerate}
\item Assume $t=2$, for each vector $(n_1,n_2)$ with $n_1+n_2=n$ with $0<n_1<n$ there is a unique matrix which is $\mathcal M = \left(n_1,0\right)$.
\item Assume $t=3$, for each vector $(n_1,n_2,n_3)$ there are $n_2-1$ matrices of type $\mathcal M =  \left(n_1, a, 0\right)$ where $0<a<n_2$. Hence,
    $$CGVN(n,3)=1(n-3)+2(n-4)+3(n-5)+\dots+(n-3)1$$
    which coincides with the given expression.
\item If we remove the null voters from a complete simple game with $t$ types of voters (having at least one null voter), we
      obtain a complete simple game with $t-1$ types of voters without null voters. Thus we have $CGVN(n,4)=$
      $$
        \sum_{i=1}^{n} CGV(n-i,3)-CGVN(n-i,3)=\sum_{i=1}^{n-4}CGV(n-i,3)-CGVN(n-i,3).
      $$
Inserting the previously obtained results
$$ CGV(n,3)=F(n+7)-\dfrac 12(n^3+2n^2+13n+26)$$
$$ CGVN(n,3)=\dfrac{(n-1)(n-2)(n-3)}6$$
for $n\ge 5$ yields $CGVN(n,4)=$
\begin{eqnarray*}
 & \sum\limits_{k=4}^{n-1} \left[ F(k+7) - \dfrac 12 (k^3+2k^2+13k+26) - \dfrac 16 (k-1)(k-2)(k-3) \right]\\
 =& F(n+8)-\dfrac 16(n^4-2n^3+26n^2+47n+132).
\end{eqnarray*}\hfill $\square$
\end{enumerate}

Starting with $n=5$ the first ten numbers of the sequence $CGVN(n,4)$ are:
$$ 1, 8, 35, 113, 303, 717, 1552, 3145, 6062, 11242.$$

Asymptotically we have
$$ \lim_{n \rightarrow \infty} \dfrac{CGVN(n,4)}{CG(n,2)} =\lim_{n \rightarrow \infty} \dfrac{F(n+8)+O(n^4)}{F(n+6)+O(n^2)}
= \lim_{n \rightarrow \infty} \dfrac{F(n+8)}{F(n+6)} = \left(\dfrac{1 + \sqrt{5}}{2}\right)^2 $$

Thus, when $n$ is fix and high enough, the number of complete games with four types of voters having vetoers and nulls (or one of the
combinations passers/nulls, vetoers/semi-vetoers, or passers/semi-passer) is almost equal to the number of complete games with two
types of voters multiplied by $\left(\dfrac{1+\sqrt{5}}{2}\right)^2$.

\section{Enumerations with more than two distinguished type of voters in a complete simple game}
\label{sec_three_distinguished_types}

Continuing the considerations from the previous two sections we study complete simple games
containing at least three of the six distinguished types from Definition~\ref{D:4types}. 
Complete simple games containing a dictator and at least another distinguished type of
voters are completely treated at the beginning of Section~\ref{sec_two_distinguished_types}.
So in the following we assume that no dictator is present, including the combination of a passer
and a vetoer, which then would be a dictator. Summarizing the results from
Section~\ref{sec_two_distinguished_types} we state that only the following seven
combinations of two distinguished type of voters are possible:
\begin{itemize}
 \item semi-passer and semi-vetoer
 \item vetoer and semi-vetoer
 \item passer and semi-passer
 \item passer and semi-vetoer
 \item vetoer and semi-passer
 \item vetoer and null voter
 \item passer and null voter
\end{itemize}
If a complete simple game contains at least one semi-passer and at least one semi-vetoer then no other distinguished
types of voters can occur. If we represent the possible combinations of two distinguished types of voters by an edge, we obtain
a quadrangle on the set $\{\text{vetoer}, \text{passer},\text{semi-vetoer},\text{semi-passer}\}$ of vertices. Thus 
there are no complete simple games with at least three distinguished types of voters.

%
%
%

\section{Conclusion and future research}
\label{sec_conclusion}

In this paper we have studied complete simple games containing at least one voter of a list of distinguished types of
voters and provided several bijections between the corresponding classes of voting systems. This contributes to the
program of numerical characterization and classification of voting systems initiated by von Neumann and Morgenstern \cite{vNM44}.

{\sloppy
It turned out that all of the counts $CGV(n,3)$, $CGN(n,3)$, $CGP(n,3)$, $CGSV(n,3)$, $CGSP(n,3)$, $CGVN(n,4)$, $CGPN(n,4)$,
$CGPSP(n,4)$, and $CGVSV(n,4)$ 
of complete simple games with distinguished and two additional types of voters
belong to the class $\Theta\left(\left(\dfrac{1+\sqrt{5}}{2} \right)^n \right)$ as well as $CG(n,2)$.} Thus, the addition
of just one of the following set of voters: vetoers, passers, nulls, vetoers and nulls, passers or nulls, semi-vetoers or
semi-passers in a complete game with two types of voters does not alter the asymptotic behavior. For a fix $n$ the
number of these games is equal to $k \cdot F(n+6) + P(n)$ where for each case $k$ is a positive constant that takes one of
the two values $\dfrac{1+\sqrt{5}}{2}$ or $\left( \dfrac{1+\sqrt{5}}{2} \right)^2$, and $P$ is a polynomial.

The exact enumeration formulas are mainly based on the previously
determined enumeration of complete simple games with two types of voters, i.e.\ $CG(n,2)$. A quite natural next step
now is to consider complete simple games with three types of voters. So far we were only able to compute the first
few exact values. For $4\le n\le 21$ the sequence $CSG(n,3)$ is given by
\begin{eqnarray*}
  && 6, 50, 262, 1114, 4278, 15769, 58147, 221089, 886411, 3806475, 17681979,\\
  && 89337562, 492188528, 2959459154, 19424078142, 139141985438, \\
  && 1087614361775, 9274721292503.
\end{eqnarray*}
For $t\ge 4$ and $n\ge 10$ we could only compute the additional values $CG(10,4)=4570902$, $CG(11,4=59776637$,
$CG(12,4)=1047858496$, $CG(13,4)=26000281487$, $CG(10,5)=412734188$, $CG(11,5)=29086472429$, and $CG(10,6)=42427707348$.

Significant sub-classes of the studied structures are those which are weighted games and, or more generally, roughly weighted
games being complete. Even though that there exist combinatorial characterizations of weighted games (see~\cite{FrMo09} and ~\cite{TaZw92})
and of roughly weighted games (see~\cite{GoSl11}), very little is known on closed enumeration formulas for these much more
restrictive voting structures, besides May's Theorem for voting systems with only one type of voters.

Using the characterization of weighted voting games with two types of voters given in \cite{FMR12ANOR} we can at least compute some
exact values of $WG(n,2)$ without generating the entire class of corresponding complete simple games with two types of voters. We would
like to remark that it took only 19~days of computation time to compute the values of $WG(n,2)$ for all $n\le 200$. Two examples
are given by $WG(100,2)=27970501$ and $WG(200,2)=851946591$. Having these numerical data at hand we observe that 
$WG(n,2)\approx0.002531n^5+O(n^4)$. And indeed it is not too hard to come up with an upper bound of $WG(n,2)\le \frac{n^5}{15}+4n^4$,
see \cite{FK12}, and a similar lower bound. The determination of an exact enumeration formula for $WG(n,2)$ seems to be an
interesting but resolvable research problem.

\bibliographystyle{elsarticle_harv}


\end{document}